\documentclass[conference]{IEEEtran}
\IEEEoverridecommandlockouts
 \usepackage[T1]{fontenc}
\usepackage[ruled,vlined]{algorithm2e}
\usepackage{cite}
\usepackage{amsmath,amssymb,amsfonts}
\usepackage{bm}
\usepackage{bbm}
\usepackage{subcaption}
\usepackage{optidef}
\usepackage{graphicx}
\usepackage{textcomp}
\usepackage{xcolor}
\usepackage[normalem]{ulem}
\usepackage{fancyhdr}
\usepackage{etoolbox}
\usepackage{booktabs}
\usepackage{hyperref}
\usepackage{tikz}
\usepackage[left=0.75in, right=0.75in, top=1in, bottom=1in]{geometry}
\usetikzlibrary{shapes, positioning}
\usepackage[tableposition=bottom]{caption}

\newtheorem{definition}{Definition}

\def\BibTeX{{\rm B\kern-.05em{\sc i\kern-.025em b}\kern-.08em
    T\kern-.1667em\lower.7ex\hbox{E}\kern-.125emX}}
\begin{document}

\IEEEpubid{\makebox[\columnwidth]{\footnotesize
\textsuperscript{~~~*}These authors contributed equally to this work.
\hfill} \hspace{\columnsep} \makebox[\columnwidth]{}}

\title{Optimal Micro-Transit Zoning via Clique Generation and Integer Programming\\
}

\author{
\IEEEauthorblockN{Hins Hu\textsuperscript{*}}
\IEEEauthorblockA{\textit{Systems Engineering} \\
\textit{Cornell University}\\
Ithaca, New York, USA \\
zh223@cornell.edu}
\and
\IEEEauthorblockN{Rhea Goswami\textsuperscript{*}}
\IEEEauthorblockA{\textit{Computer Science} \\
\textit{Cornell University}\\
Ithaca, New York, USA \\
rkg62@cornell.edu}
\and
\IEEEauthorblockN{Hongyi Jiang}
\IEEEauthorblockA{\textit{Systems Engineering} \\
\textit{City University of Hong Kong}\\
Kowloon, Hong Kong, China\\
hongyi.jiang@cityu.edu.hk}
\and
\IEEEauthorblockN{Samitha Samaranayake}
\IEEEauthorblockA{\textit{Civil and Environmental Engineering} \\
\textit{Cornell University}\\
Ithaca, New York, USA\\
samitha@cornell.edu}
}

\maketitle

\begin{abstract}
Micro-transit services offer a promising solution to enhance urban mobility and access, particularly by complementing existing public transit. However, effectively designing these services requires determining optimal service zones for these on-demand shuttles, a complex challenge often constrained by operating budgets and transit agency priorities. This paper presents a novel two-phase algorithmic framework for designing optimal micro-transit service zones based on the objective of maximizing served demand. A key innovation is our adaptation of the shareability graph concept from its traditional use in dynamic trip assignment to the distinct challenge of static spatial zoning. We redefine shareability by considering geographical proximity within a specified diameter constraint, rather than trip characteristics. In Phase 1, the framework employs a highly scalable algorithm to generate a comprehensive set of candidate zones. In Phase 2, it formulates the selection of a specified number of zones as a Weighted Maximum Coverage Problem, which can be efficiently solved by an integer programming solver. Evaluations on real-world data from Chattanooga, TN, and synthetic datasets show that our framework outperforms a baseline algorithm, serving 27.03\% more demand in practice and up to 49.5\% more demand in synthetic settings.
\end{abstract}


\section{Introduction}
Public transit systems are fundamental to sustainable urban development, contributing to reduced emissions and energy consumption \cite{b1} while fostering economic development and social equity \cite{b2}. However, inadequate public transit investments along with car-centric planning have resulted in very restricted transit coverage in most U.S. cities, leaving many communities with limited access to essential services and opportunities via public transit. \emph{Micro-transit}, a small-scale demand-responsive transit service characterized by flexible routes and adaptable scheduling, has recently emerged as a promising complement to traditional fixed-route mass transit---with the potential for improving access in areas with restricted coverage, such as lower-density areas \cite{b3}.

One of the primary decisions for an operator when planning a micro-transit service is to decide where to operate the service given budgetary constraints and other operational requirements/constraints. This optimal micro-transit zone design problem is a complex optimization problem and typically solved either directly via expert opinion or through heuristic methods. While zoning methodologies exist, many are adapted from different contexts like freight logistics  (e.g., using clustering techniques) \cite{b9} or, when specific to micro-transit, may prioritize objectives like equity without directly optimizing for demand coverage based on origin-destination patterns. Consequently, these approaches fall short in identifying zone configurations that maximize service coverage under operator constraints.

To address the zoning problem, this paper proposes a novel, optimal, and highly scalable algorithmic framework. The core idea of this framework is inspired by the concept of the shareability graph, originally developed for taxi-sharing systems \cite{b4}. The original shareability graph is a spatiotemporal network abstraction used to characterize which taxi requests can share a taxi given user constraints such as waiting time, that can be used to efficient compute optimal sharing strategies on large-scale taxi trip datasets. It has subsequently been applied to other operations research problems in transportation systems, such as high-capacity ride-pooling assignment \cite{b12} and school bus routing \cite{b5}

In this paper, we adapt the shareability graph concept specifically for the micro-transit zoning problem. Instead of focusing on whether specific trips can be shared (as in the original concept), our key modification is to focus on whether locations are close enough to be grouped. In our model, we represent locations as ``nodes'', where each node stands for a small neighborhood and captures the travel demand originating or ending in that area. Two nodes are considered ``shareable" if the shortest-path distance between them is within a specified diameter constraint $D$, which defines the maximum geographical extent permitted for any single service zone.

Based on this modified shareability concept, we design a two-phase algorithmic framework to find optimal zones. Phase 1 is an algorithm called \textbf{\textsc{CliqueGen}}, which leverages the modified shareability structure to iteratively construct all valid candidate zones, represented as ``cliques'' of mutually shareable nodes, that adhere to the diameter constraint $D$. If a zone is selected, all travel demand where both origin and destination fall within that zone's constituent nodes will be served. Phase 2 formulates the zone selection task as a Weighted Maximum Coverage Problem, taking the candidate zones from Phase 1 as input.  An Integer Linear Program (ILP) called \textbf{\textsc{ZoningILP}} is then solved to select the $m$ zones that collectively cover the maximum possible intra-zone travel demand. The effectiveness of this framework was validated through numerical experiments using real-world trips from Chattanooga, TN, USA, demonstrating practical applicability in an operational setting. We complement this with synthetically generated datasets that enable systematic sensitivity analysis across varying problem parameters, including network size and diameter constraints, which would be difficult to control with real-world data alone.

\medskip
work for designing effective micro-transit zones, contributing
novel modeling techniques, an efficient solution approach, and its validation. Several key points are summarized as follows:
\begin{enumerate}
    \item We propose a highly scalable two-phase algorithm to optimally solve the micro-transit zoning problem.

    \item We validate the algorithm through numerical experiments, demonstrating a substantial increase in total demand served compared to a baseline algorithm.
\end{enumerate}

\section{Related Work}

This section first examines existing methodologies for transportation zoning, with a particular focus on approaches relevant to micro-transit design. Second, it discusses the development and application of shareability graphs within transportation research, a concept central to our proposed framework.

\textit{a) Transportation Zoning Approaches:} Zoning is fundamental to urban planning and infrastructure management, influencing land use patterns \cite{b6} and correlating with urban development trajectories \cite{b7}.  However, poorly designed or outdated zoning can exacerbate urban inequalities \cite{b8}. Within transportation planning, computational approaches have been developed to automate or optimize zone design. For instance, in freight transportation planning, methodologies have been proposed to update zoning systems by clustering geographically similar spatial units using techniques like K-means clustering \cite{b9} or multi-objective genetic algorithms for improved scalability \cite{b10}. In micro-transit contexts, \cite{erdogan2024} developed a framework incorporating expert judgment for zone evaluation, while \cite{ng2024} focused on operational dispatching within pre-determined zones using reinforcement learning.

While the benefits of micro-transit are increasingly recognized, the specific problem of computationally designing optimal micro-transit zones is less explored. Existing research often focuses on different objectives or employs distinct methodologies. For example, Bonner \& Miller-Hooks \cite{b11} developed a mixed-integer programming approach to generate equitable micro-transit zones, maximizing equity-based objective functions that prioritize fair service distribution across populations and geographic areas, rather than directly targeting demand coverage. Our work differs by introducing a optimal and highly-scalable graph-based methodology built upon an efficient candidate zone generation and selection process. 

\textit{b) Shareability Graph:} The concept of the shareability graph, introduced by Santi et al. \cite{b4} for taxi sharing, has proven influential in addressing various transportation optimization problems. This graph-abstraction captures the potential for individual requests to be grouped together based on spatio-temporal constraints. It was initially used to quantify the benefits of taxi sharing by efficiently computing optimal sharing strategies \cite{b4}. Subsequently, variations of the shareability graph have been applied to complex operational problems in many specific contexts, including, most prominently, the real-time matching of riders and drivers in high-capacity ride-pooling systems \cite{b12}.

Guo et al. \cite{b5} adopted the shareability graph to solve large-scale school bus routing problems via a decomposition based on feasible bus trips. Shah et al. \cite{Shah2019} extended shareability graph approaches with neural network-based value functions to overcome the myopic nature of ride-pooling assignments. Vazifeh et al. \cite{Vazifeh2018} adapted the shareability graph concept to create ``vehicle-sharing networks" that model the sharing of vehicles rather than rides. Relevant to micro-transit's role in complementing mass transit, Edirimanna et al. \cite{b13} modified the graph construction to model ride-sharing for first-and-last-mile trips in a multi-modal setting, defining feasibility based on subgraph structures and quality-of-service constraints. 

While such work advances the operational planning of ride-sharing services, our research addresses the distinct, foundational challenge of designing the geographical service zones themselves. We apply the shareability graph concept differently, shifting to a static, node-level perspective based on geographical proximity. This novel adaptation allows us to specifically tackle the a priori zoning problem, providing a method to strategically define optimal service areas to allow effective operational deployment.

\section{Problem Description} \label{sec:problem}
We formulate the \textbf{micro-transit zoning problem} as a combinatorial optimization problem with the following inputs:
\begin{enumerate}
    \item A planar directed graph $G = (V, E)$ abstracted from an urban region, where each node $v \in V$ represents a small neighborhood that aggregates all travel demand from and to that area. An edge $(i, j) \in E$ exists only if node $i$ and $j$ are connected by a segment of a road. The edge weight $w(i, j)$ is defined as the center-to-center distance between node $i$ and $j$.
    \item A demand table $\bm{d}: V^2 \rightarrow \mathbb{R}_+$, where each entry $d(i, j)$ represents a nominal travel demand from node $i$ to $j$.
    \item A distance matrix $\bm{c}: V^2 \rightarrow \mathbb{R}_+$, where each entry $c(i, j)$ denotes the shortest-path distance from node $i$ to $j$ in graph $G$. 
    \item A maximum diameter $D \in \mathbb{R}_+$ constrains the spatial extent of each zone, defined as the maximum pairwise distance between any two nodes assigned to the zone.
    \item An integer $m \in \mathbb{N}$ indicating the number of zones to be set up.
\end{enumerate}
The goal is to select $m$ micro-transit zones, where each zone is defined as a subset of connected nodes $S \subseteq V$ satisfying the diameter constraint $\max_{i, j \in S} c(i, j) \leq D$, such that the total intra-zone demand served in Equation \ref{eq:obj} is maximized.
\begin{equation} \label{eq:obj}
\sum\limits_{i,j \in V} d(i, j) \; \mathbb{I} (i, j \text{ in the same zone})  
\end{equation}

\section{Methodology}
\begin{figure}[!htbp]
    \centering
    \includegraphics[width=0.8\linewidth]{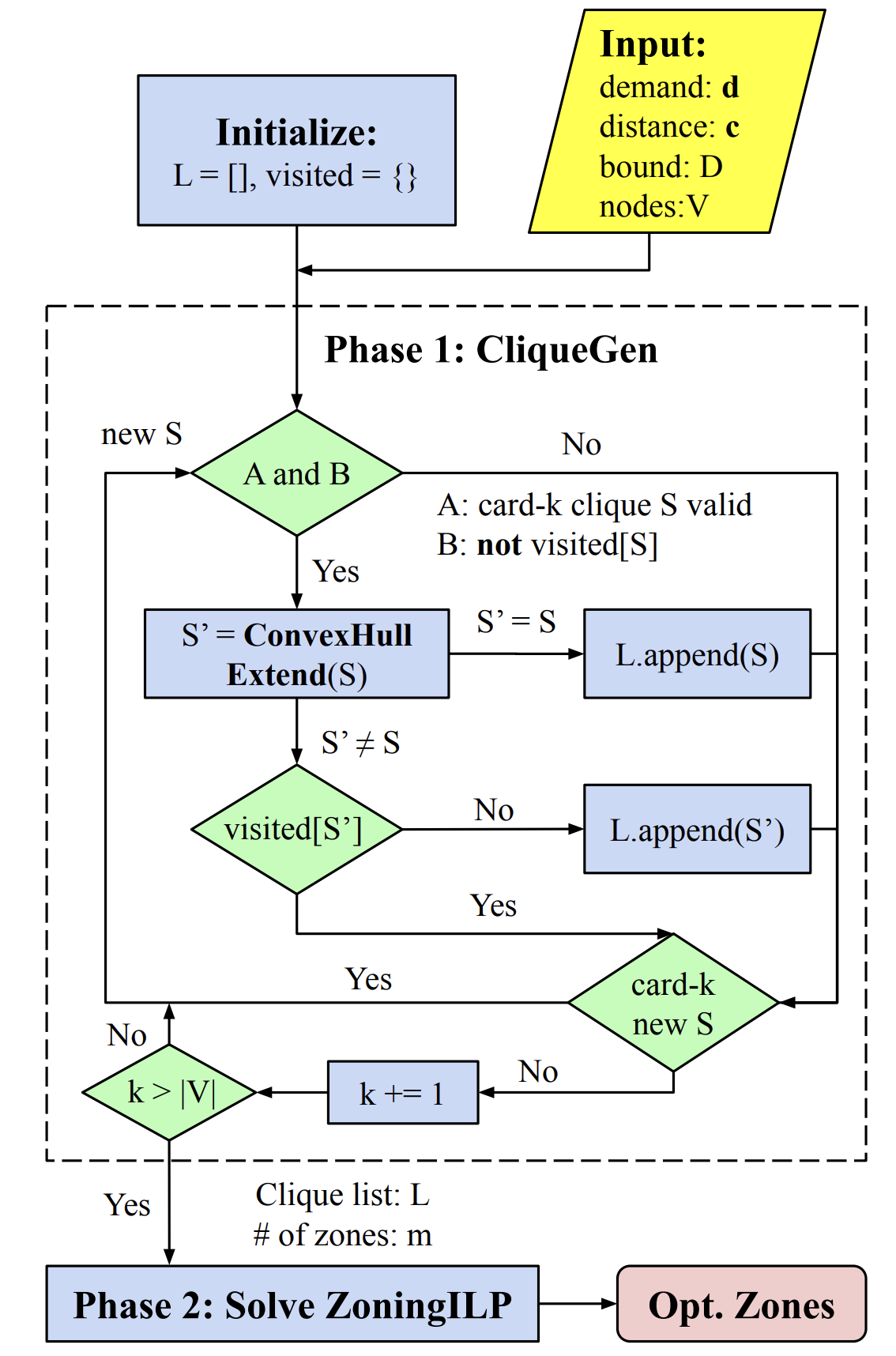}
    \caption{A schematic diagram of the two-phase algorithmic framework.}
    \label{fig:diagram}
\end{figure}

We approach the micro-transit problem via a two-phase algorithmic framework, which is illustrated in the schematic diagram in Figure \ref{fig:diagram}.

In Phase 1, we design a candidate zone generation algorithm, called \textsc{CliqueGen}, to efficiently enumerate all valid zones with diameters less than $D$, with the full procedure described in Section \ref{sec:clique_gen}. Inspired by the definition of shareability graph \cite{b4}, the algorithm iteratively constructs cliques\footnote{Following prior work on shareability graphs, we adopt the term \textit{clique} for our subgraph structure of interest, though our definition differs from the classical graph-theoretic one that a clique is a complete subgraph.} -- special subgraph structures -- of cardinality $k$ by progressively extending cliques of cardinality $k-1$. A clique can be regarded as a valid zone and is defined below.
\vspace{0.1cm}
\begin{definition}[$D$-bounded Clique] \label{def:clique}
Given a graph $G = (V, E)$, a $D$-bounded clique is a subset of nodes $S \subseteq V$, such that the longest shortest-path distance between any two nodes $i, j \in S$, called the diameter, is less than $D$. The cardinality of a clique is the number of nodes in it, namely $|S|$.
\end{definition}
\vspace{0.1cm}

Although the number of valid zones grows exponentially when the diameter constraint is relaxed, \textsc{CliqueGen} leverages a pruning sub-routine called \textsc{CovexHullExtend} to exclude a large number of zones that do not make physical sense and are dominated by other zones that we keep. Any such dominated zone does not cover all nodes within its convex hull, and thus cannot contribute to an optimal solution to our problem. This significantly reduces the computational overhead and enables practical tractability, as demonstrated in Section \ref{sec:res}. The number of candidate zones generated by \textsc{CliqueGen} for different parameters is detailed in Table \ref{tab:scalability}. 

In Phase 2, we formulate an Integer Linear Program (ILP), denoted by $\textsc{ZoningILP}$, to select the optimal $m$ zones given the full list of cliques generated in Phase 1. Such an optimization problem is a variant of the classical Weighted Maximum Coverage Problem (WMCP) \cite{b14}. The full formulation and discussions are in Section \ref{sec:ilp}.

\subsection{Candidate Zone Generation} \label{sec:clique_gen}
Given a diameter bound $D$, \textsc{CliqueGen} builds up the list of cliques in an iterative manner. It begins by initializing a list $L$ with all cardinality-$1$ cliques (i.e., single nodes). Then, at any arbitrary iteration $k$, \textsc{CliqueGen} attempts to construct all possible cliques of cardinality $k$ that can be directly extended from at least one cardinality-($k$-1) clique. Specifically, if a clique such as $\{1, 5, 9\}$ is deemed invalid at iteration $k-1$, any superset, like $\{1, 5, 9, v\}$ with $v$ being any other nodes, will be excluded by the algorithm for further consideration. A new cardinality-$k$ clique formed by adding node $v$ to an existing clique $S \in V^{k-1}$ is retained in list $L$ only if $\max\{c(i, v), \; c(v, i)\} \leq D$ for all $i \in S$. This condition ensures that $v$ lies within the diameter constraint relative to all members of $S$, effectively enforcing a local geographic proximity requirement.

This idea is analogous to the the notion of shareability in prior works \cite{b4} \cite{b12} \cite{b13}, where trip-level shareability is verified based on the temporal or spatial proximity between origin-destination pairs. However, our setting is even simpler in two key ways: (1) \textsc{CliqueGen} assesses the shareability at the node level rather than the trip level, and (2) pairwise node-level shareability suffices to guarantee clique-level validity. The latter is in contrast to trip-based models, where verifying group shareability typically requires solving an NP-hard Pick-up and Delivery Problem (PDP) for all origins and destinations, whereas our clique validation takes only at most $\binom{|S|}{2}$ distance checks.

At an iteration $k$, once \textsc{CliqueGen} identifies a new clique $S^\prime = S \, \cup \{v\}$ from an existing clique $S$, it invokes a sub-routine called \textsc{ConvexHullExtend} to check whether the new clique $S^\prime$ can be further extended by adding nodes within its convex hull, $\text{Conv}(S^\prime)$, geographically. If \textsc{ConvexHullExtend} can find a set of additional nodes $Q$ contained in $\text{Conv}(S^\prime)$, the convex-hull-extended clique $S^\prime \cup Q$ is added to the list $L$, and clique $S^\prime$ is discarded. Otherwise, clique $S^\prime$ is added to $L$. Once all cardinality-$k$ clique extensions are processed, \textsc{CliqueGen} moves to iteration $k+1$. The algorithm terminates when $k$ cannot increase anymore, which occurs when cliques are over-extended (i.e., new nodes are far away) and the diameter constraint is near violation. The pseudo-code for \textsc{CliqueGen} is provided in Algorithm \ref{algo:clique_gen}.
\begin{figure}[htbp]
    \centering
    \begin{tikzpicture}[scale=1.5]

    \coordinate (N1) at (0, 0);
    \coordinate (N2) at (1.5, 0);
    \coordinate (N3) at (2, 1.5);
    \coordinate (N4) at (0.5, 2);
    \coordinate (N5) at (-0.3, 1);

    \coordinate (N6) at (0.5, 0.6);
    \coordinate (N7) at (1.3, 0.5);
    \coordinate (N8) at (0.7, 1.4);

    \coordinate (N9) at (-0.3, 0.3);
    \coordinate (N10) at (2.5, 1);
    \coordinate (N11) at (1.5, 1.9);

    \draw[draw=blue, thick] (N1) -- (N2) -- (N3) -- (N4) -- (N5) -- cycle;

    \foreach \pt/\name in {N1/1, N2/2, N3/3, N4/4, N5/5, N8/8} {
        \fill[blue] (\pt) circle (1.2pt);
        \node[black, font=\small, above left] at (\pt) {\name};
    }
    
    \foreach \pt/\name in {N6/6, N7/7} {
        \fill[red] (\pt) circle (1.2pt);
        \node[black, font=\small, above left] at (\pt) {\name};
    }
    
    \foreach \pt/\name in {N9/9, N10/10, N11/11} {
        \fill[green!70!black] (\pt) circle (1.2pt);
        \node[black, font=\small, above left] at (\pt) {\name};
    }
    \end{tikzpicture}
    \caption{An example: clique $S_1 = \{1, 2, 3, 4, 5, 8\}$ are valid zones but dominated by clique $S_2 = \{1, 2, 3, 4, 5, 6, 7, 8\}$. Integrating nodes $9$, $10$, or $11$ to $S_1$ requires shareability checks, but $6$ and $7$ should be automatically included through \textsc{ConvexHullExtend}.}
    \label{fig:convex_hull}
\end{figure}
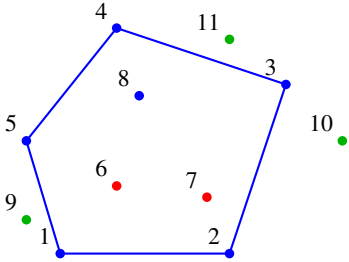

The rationale behind the convex hull extension is that a ``hollow'' clique, one that omits nodes contained within its convex hull, will always be dominated by the zone that includes the omitted nodes, as illustrated in Figure \ref{fig:convex_hull}. This is because it has the same diameter as the ``hollow'' clique (i.e., it is also feasible), but potentially serves more travel demand. Therefore, \textsc{CliqueGen} will only include \textit{convex} cliques, formalized in Definition \ref{def:convex_clique}, into the list $L$.
\begin{algorithm}[!htbp]
\caption{CliqueGen}
\label{algo:clique_gen}
\KwIn{Nodes $V$, Distances $\bm{c}$, Max Diameter $D$}
\KwOut{Clique List $L$}

$L \gets [\{v\} \mid v \in V]$

$k \gets 1$

\While{true}{
    \If{no valid clique of cardinality $k$}{
        \textbf{break}
    }
    \ForEach{clique $S$ of cardinality $k$ in $L$}{
        \ForEach{$v \in V \setminus S$}{
            \If{$\max\{\bm{c}[v,u], \bm{c}[u,v]\} \leq D, \; \forall u \in S$}{
                $S' \gets \textsc{ConvexHullExtend}(S, V, \bm{c}, D)$ \\ 
                \eIf{$S' = \varnothing$}{
                    $L.\text{append}(S)$
                }{
                    $L.\text{append}(S')$
                }
            }
        }
    }
    $k \gets k + 1$
}
\Return $L$

\end{algorithm}
\begin{definition}[Convex Clique] \label{def:convex_clique}
A clique $S$ is convex if it is \textit{closed} under the convex hull extension in the node set $V$. Equivalently, $S = \{v \in \text{Conv}(S) \, | \, v \in V\}$.
\end{definition}

The sub-routine \textsc{ConvexHullExtend} is executed as follows: given a clique $S$ and the node set $V$, it iteratively examines each node $v \in V \setminus S$ and determines if $v$ is shareable with all nodes in $S$, using the same pairwise shareability check as in \textsc{CliqueGen}. Since the internal shareability of $S$ is already validated, each new candidate node $v$ only requires $O(|S|)$ checks to complete the validation. If $v$ is shareable with the entire $S$, an oracle is further invoked to check whether $v$ is contained in $Conv(S)$ before integrating $v$ into $S$. The original clique $S$ keeps expanding along the loop until no more new node is found within $\text{Conv}(S)$. Given that the ``point-in-convex-hull'' oracle takes $O(|S|)$ steps, the overall time complexity of \textsc{ConvexHullExtend} is $O(|S|^2)$. The pseudo-code is provided in Algorithm \ref{algo:convex_hull}.
\begin{algorithm}[!htbp]
\caption{ConvexHullExtend}
\label{algo:convex_hull}
\KwIn{Clique $S$, Nodes $V$, Distances $\bm{c}$, Max Diameter $D$}
\KwOut{Convex-hull-extended clique $S'$}

\While{$V$ is not empty}{
    $v \gets V.\text{pop}()$\\ 
    \If{$\max\{\bm{c}[v,u], \bm{c}[u,v]\} \leq D$ for all $u \in S$}{
        \If{$v \in \text{Conv}(S)$}{
            $S \gets S \cup \{v\}$
        }
    }
}
$S' \gets S$

\Return $S'$
\end{algorithm}

It is important to note that a naive implementation of \textsc{CliqueGen} can be highly inefficient due to extensive re-computation in the following aspect: the algorithm may repeatedly examine the same set of nodes as it checks the formation of an identical valid clique from different lower-cardinality cliques. To address this issue, we design a data structure called \texttt{visited} that maps the sorted tuple of cliques to a binary flag indicating whether the combination has been explored in previous steps. \textsc{CliqueGen} always checks \texttt{visited} before extending an existing clique $S$ with a new node $v$, or with a new set of nodes $Q$ contained in $\text{Conv}(S)$. If \texttt{visited} returns true for a combination, the step is skipped to save computation. Similarly, whenever \textsc{CliqueGen} inserts a newly found clique $S$ into the list $L$, it updates the map with \texttt{visited[sorted($S$)] = 1}.

\subsection{Optimal Zoning via Maximum Coverage} \label{sec:ilp}
Given the comprehensive list $L$ of all the candidate valid zones  generated by \textsc{CliqueGen}, the goal in Phase 2 is to select $m$ zones from $L$ that maximizes the total demand served, as described in Equation \ref{eq:obj}. Since the demand table is defined over node pairs, allowing zones to overlap can lead to higher total coverage and thus better solutions. This problem reduces to a Weighted Maximum Coverage Problem with $m$ sets. 

Let $x_S \in \{0, 1\}$ be a binary decision variable indicating whether a candidate valid zone $S \in L$ is selected. Let $y_{ij} \in \{0, 1\}$ be a binary decision variable indicating whether both nodes $i, j \in V$ are jointly covered by at least one selected zones. For each zone $S \in L$, define a binary vector $z_{S} = [z_{S}^1, \dots, z_{S}^i, \dots, z_{S}^{|V|}] \in \{0, 1\}^V$ capturing the mapping between zone $S$ and all the nodes in $V$. In particular, $z_S^i = 1$ if node $i \in S$, and $0$ otherwise. These vectors $\{z_{S}\}_{S \in L}$ are provided as input as they can be directly derived from the clique list $L$ generated in Phase 1. With the notations properly defined, we formulate \textsc{ZoningILP}, consisting of the objective function and constraints from (2a) to (2e), to solve this maximum coverage problem exactly. 
\begin{maxi!}|s|[2]<b>
    {\bm{x}, \bm{y}}
    {\sum_{i, j \in V} d_{ij} \, y_{ij}}{}{}
    \addConstraint{\sum_{S \in L} x_S}{\leq m}{} \label{c:master_budget}
    \addConstraint{y_{ij}}{\leq \sum_{S \in L} z_S^i \, z_S^j \, x_S, \quad}{\forall i, j \in V} \label{c:w_ij}
    \addConstraint{x_S}{\in \{0,1\}, \quad}{\forall S \in L}
    \addConstraint{y_{ij}}{\in \{0, 1\}, \quad}{\forall i, j \in V}
\end{maxi!}

The objective function in (2a) collects the travel demand of all node pairs covered by at least one zone (with $y_{i, j} = 1$). Constraint (2b) limits the number of selected zones to be $m$. The inequality is tight at optimality due to the maximization. Constraints (2c) link $\bm{x}$ and $\bm{y}$, forcing $y_{ij}$ to be $0$ if node $i$ and $j$ are not assigned simultaneously to any zone $S$ (with at least one of $z_S^i$, $z_S^j$, and $x_S$ being $0$). In the maximization setting, it also raises $y_{ij}$ to $1$ if node $i$ and $j$ are assigned to at least one zone (with $z_S^i = 1$, $z_S^j = 1$, and $x_S = 1$ for at least one $S \in L$). Constraints (2d) and (2e) enforce the integrality of decision variables. 


The size complexity of \textsc{ZoningILP} is $O(|L| + |V|^2)$ as it involves $O(|L| + |V|^2)$ decision variables and $O(|V|^2)$ constraints. While the power set of $V$ has size $O(2^{|V|})$, providing a natural upper bound on the number of candidate cliques $|L|$, this worst case is rarely realized in practice. The diameter constraint $D$, the pruning effect of \textsc{ConvexHullExtend}, and implementation techniques discussed in Section \ref{sec:clique_gen} all substantially reduce $|L|$, making the ILP tractable for medium-to-large size city, as demonstrated in Section \ref{sec:res}.

Moreover, we conjecture that, due to the convexity property of valid zones in Definition \ref{def:convex_clique}, the theoretical upper bound on $L$ could be \textit{sub-exponential} or even \textit{polynomial} in sufficiently regular graphs (e.g., uniform grids). A formal complexity analysis for graphs with additional structure is beyond the scope of this work and is left for future research.  

\section{Numerical Experiments} \label{sec:res}
We conduct extensive numerical experiments using both real-world travel demand data from Chattanooga, TN and synthetic datasets based on a pruned Delaunay graph.

To benchmark the performance of our proposed method, we implement a simple baseline algorithm based on a construction heuristic, referred to as \textsc{SimpleZoning}, which generates $m$ high-quality zones under the same diameter constraint $D$ and in the same problem setting.

\textsc{SimpleZoning} begins by selecting node pair $(u, v)$ with the highest travel demand as the initial zone seed $S = \{u, v\}$. It then incrementally expands $S$ by adding one node $i$ at a time that yields the greatest marginal gain in total demand served, computed as in Equation \ref{eq:demand_gain}, until no further node can be added without violating the diameter constraint $D$. Next, the algorithm removes $S$ from $V$ and repeats the same process until $m$ zones are found.
\begin{equation} \label{eq:demand_gain}
\Delta_i = \max\limits_{i \in V \setminus S} \; \sum_{v \in S} \; d[i, v] + d[v, i]
\end{equation}

\subsection{Real-World Data}
For the real-world experiments, we obtain a set of origin-destination estimates derived from a proprietary dataset. Due to the sensitivity of anonymized individual trip records, we are unable to publicly release the dataset. However, we share our \href{https://github.com/hins-hu/clique-gen}{codebase} and the synthetic datasets for reproducibility. We apply a geo-fencing filter to retain only trips with both origins and destinations within the city boundary. Additionally, we exclude very short trips whose origins and destinations are within a walkable distance (less than $500$ meters). After data cleaning and pre-processing, a total number of 33,256 trips are included in the experiments.

\begin{figure}[htbp!]
    \centering
    \includegraphics[width=\linewidth]{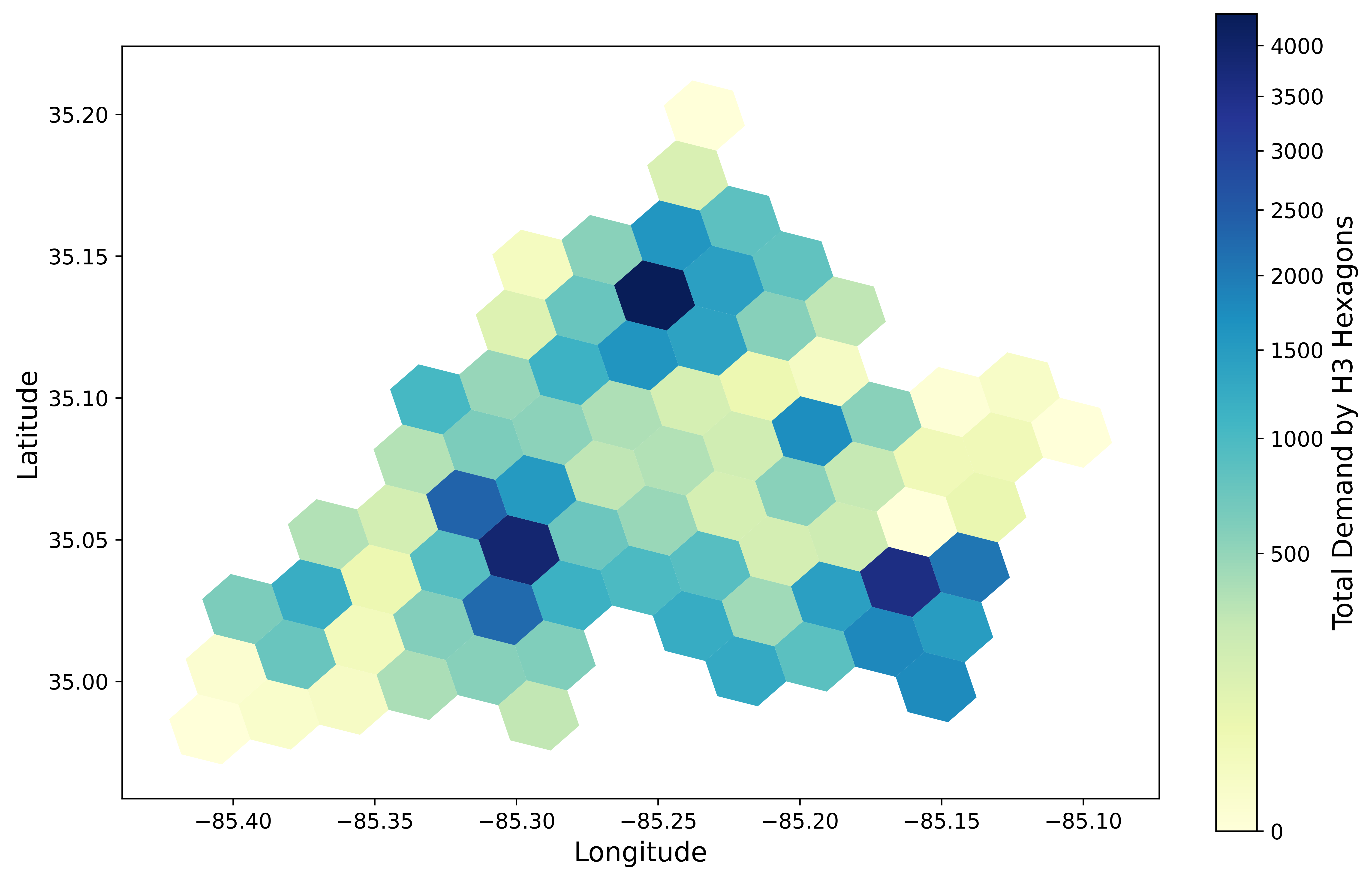}
    \caption{Visualization of H3-aggregated demand pattern in Chattanooga, TN, USA}
    \label{fig:vis_demand}
\end{figure}

To spatially aggregate the travel demand into small, contiguous neighborhoods, with each corresponding to a node as defined in Section~\ref{sec:problem}, we use the open-source \href{https://github.com/uber/h3}{Hexagonal Hierarchical Spatial Index (H3)} developed by Uber. H3 partitions geographic space into a grid of hexagons at multiple resolutions, originally designed to cover the entire planet. Since the resolution is restricted to 16 discrete levels and cannot be finely tuned between levels, we fix the resolution to a reasonable level 7, which divides the city into 78 connected hexagonal cells. Each hexagon is treated as a node, representing a small neighborhood. Trips are assigned to hexagons by mapping their origin and destination coordinates to the corresponding cells. The resulting aggregated demand pattern, where each cell represents the sum of all trips either originating from or ending at that cell, is visualized in Figure~\ref{fig:vis_demand}.

To better reflect practical conditions, we use pre-computed nominal driving times on the real-world road network for the distance matrix $\bm{c}$. Under this consideration, the maximum diameter represents the longest driving time between the central locations of any two hexagons within a zone. We set the number of zones to $m = 2$ and the maximum diameter to $D = 480$ seconds, equivalent to an 8-minute drive.

\begin{figure}[htbp!]
    \centering
    \includegraphics[width=\linewidth]{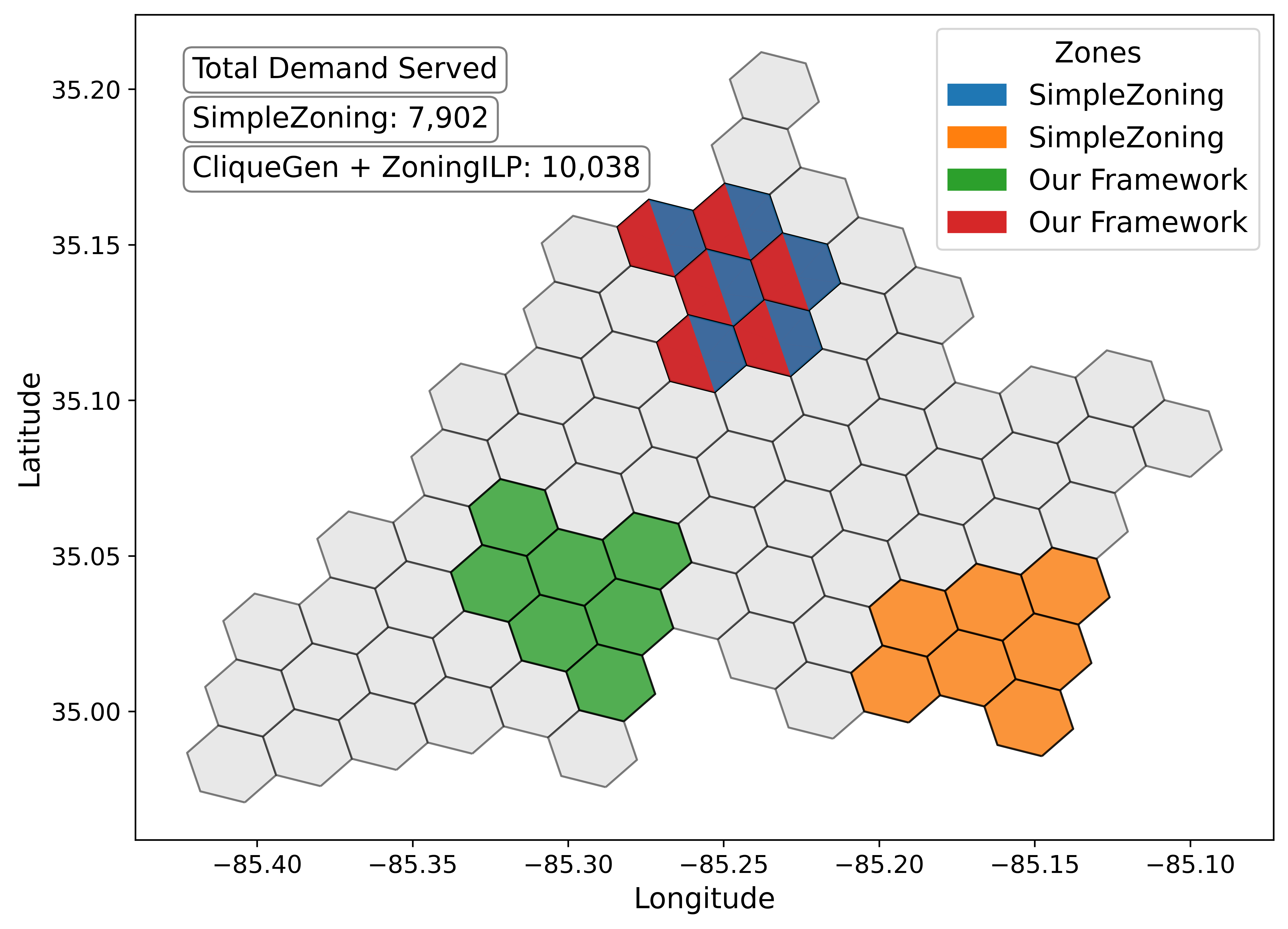}
    \caption{Micro-transit zones generated by both algorithms}
    \label{fig:real-world-zones}
\end{figure}

The zoning results of \textsc{CliqueGen} + \textsc{ZoningILP} \textbf{(our method)} and \textsc{SimpleZoning} (the baseline algorithm) are shown in Figure~\ref{fig:real-world-zones}. Both approaches select zones that align well with the demand hot spots illustrated in Figure~\ref{fig:vis_demand}. Notably, our method serves 27.03\% more trips than the baseline. Given that the budget only allows for two zones (i.e., $m = 2$), identifying the optimal subset of nodes among three prominent high-demand areas is nontrivial even for expert city planners. While the baseline converges to a local optimum by selecting nodes in the lower-right (orange) region, our optimization framework correctly prioritizes the lower-left (green) region as more beneficial for establishing a micro-transit zone.

\subsection{Sensitivity Analysis in Synthetic Datasets}
In the synthetic setting, we construct a Delaunay graph with $|V|$ nodes uniformly distributed within a $10 \times 10$ bounding box. Each undirected edge is converted into two directed edges with a probability of $0.2$, and then a pruning step removes existing edges independently with probability $0.2$. This pre-processing emulates the mixed road types and irregular connectivity patterns commonly observed in urban environments.

Compared to real-world datasets constrained by the H3 spatial index, the synthetic setting offers finer-grained control over key parameters. To systematically evaluate the performance and scalability of our framework, we vary the number of nodes $|V| \in \{50, 100, 150, 200\}$ and the maximum diameter constraint $D \in \{1.5, 2, 2.5, 3\}$. The number of zones is set to $m = 4$.

\begin{figure}[h]
    \centering
    \begin{subfigure}[b]{0.49\linewidth}
        \centering
        \includegraphics[width=\linewidth]{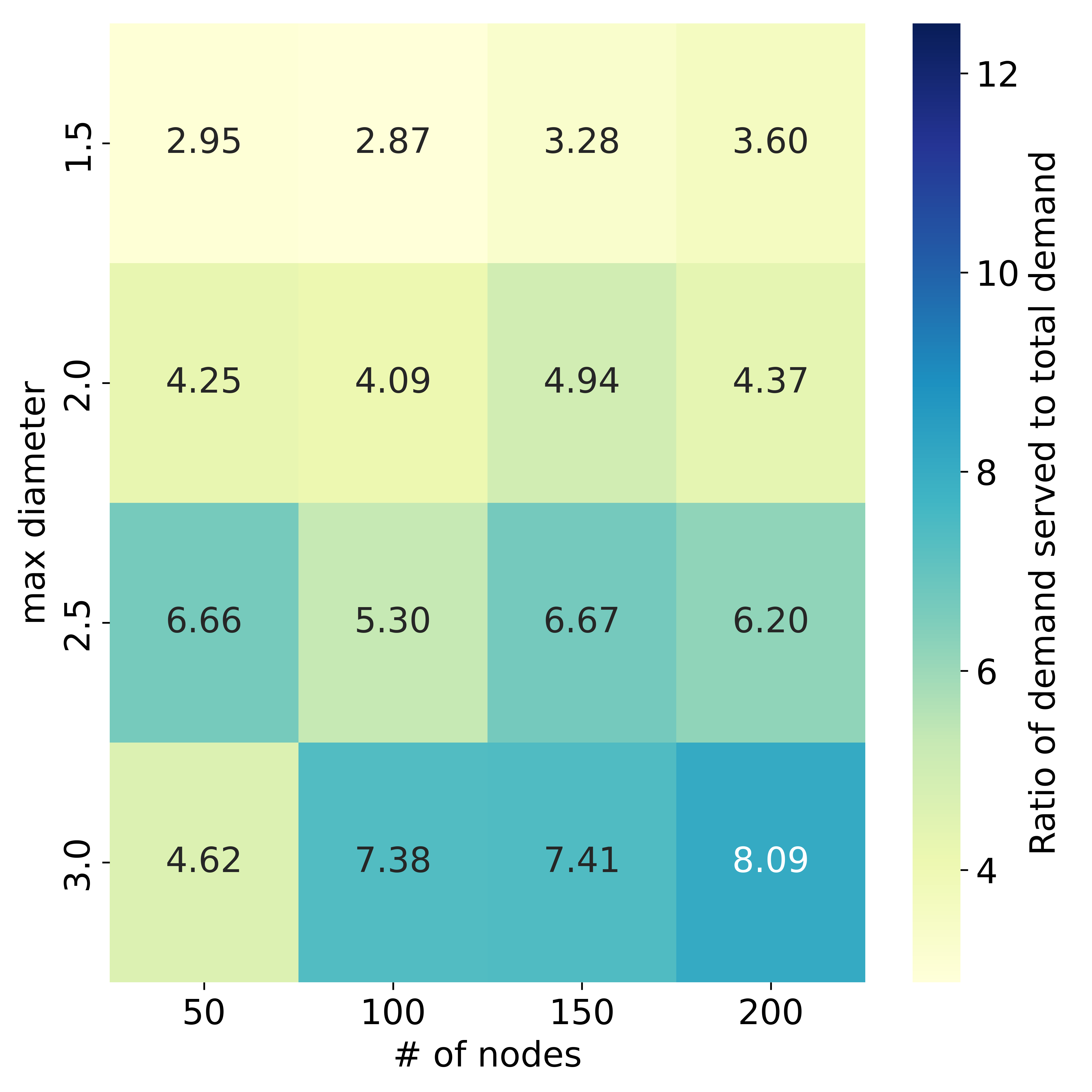}
        \caption{\textsc{SimpleZoning}}
        \label{fig:heatmap-baseline}
    \end{subfigure}
    \hfill
    \begin{subfigure}[b]{0.49\linewidth}
        \centering
        \includegraphics[width=\linewidth]{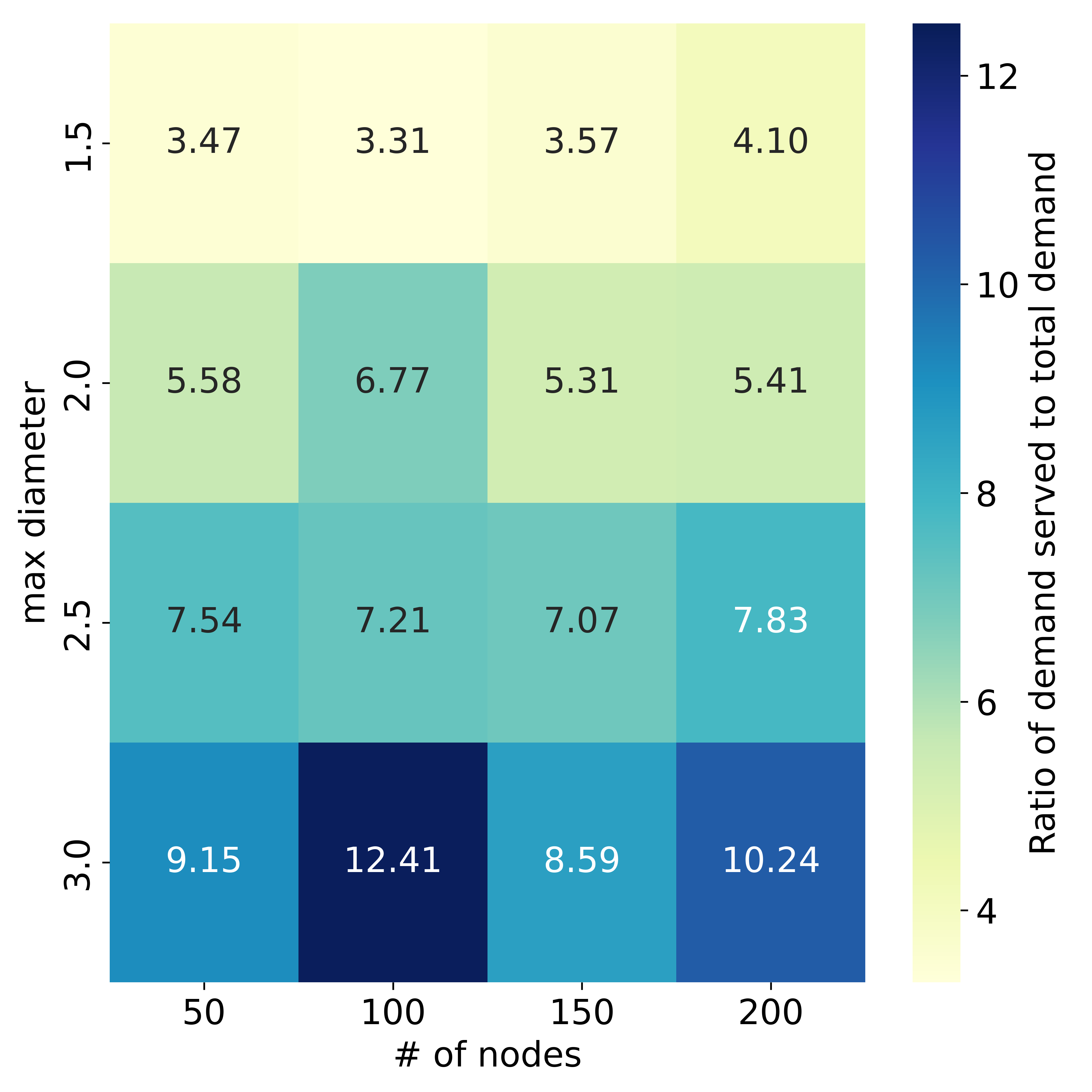}
        \caption{\textsc{CliqueGen} + \textsc{ZoningILP}}
        \label{fig:heatmap-our}
    \end{subfigure}
    \caption{Ratios of demand served to total demand with varying numbers of nodes and maximum diameters}
    \label{fig:heatmap}
\end{figure}

Figure~\ref{fig:heatmap} presents the ratio of demand served within selected zones to the total demand generated across the entire graph. Our framework consistently outperforms the baseline across all combinations of maximum diameter and node count, achieving an average improvement of 20.44\% and a maximum improvement of 49.5\% (at $D = 3$, $|V| = 50$) in total demand served. As the diameter constraint $D$ increases, the demand served by \textsc{CliqueGen} + \textsc{ZoningILP} increases monotonically for all values of $|V|$. In contrast, \textsc{SimpleZoning} is less robust as the ratio of demand served drops from $6.66\%$ to $4.62\%$ when $D$ increases from $2.5$ to $3$.

Table~\ref{tab:scalability} shows the actual number of cliques generated by \textsc{CliqueGen} across all combinations of both parameters. At the largest setting with $200$ nodes, the range can cover the entire area of a medium-to-large size city in the U.S., given that each hexagon in H3 at resolution $7$ spans approximately an $1$ km by $1$ km neighborhood. The total number of cliques in this setting expands up to $25,000$. Despite this scale, the instance remains tractable for \textsc{ZoningILP}, as the underlying Weighted Maximum Coverage Problem belongs to a relatively simple class of NP-hard problems \cite{b14}. Moreover, Figure~\ref{fig:heatmap} empirically confirms that problems of this size are solvable within reasonable computational limits. These findings support the conclusion that our framework is scalable for addressing micro-transit zoning in most medium-sized cities worldwide.
\begin{table}[!htbp]
    \centering
    \begin{tabular}{c|cccc}
        \toprule
        \textbf{Max Diameter} & \textbf{50 Nodes} & \textbf{100 Nodes} & \textbf{150 Nodes} & \textbf{200 Nodes} \\
        \midrule
        1.5 & 123 & 603 & 1314 & 3308 \\
        2.0 & 199 & 1357 & 4804 & 14182 \\
        2.5 & 293 & 3134 & 14665 & 56050 \\
        3.0 & 495 & 8257 & 44324 & 253552 \\
        \bottomrule
    \end{tabular}
    \caption{Total number of candidate zones generated by \textsc{CliqueGen} in Phase 1 with varying numbers of nodes and maximum diameters}
    \label{tab:scalability}
\end{table}

\section{Conclusion \& Future Directions}
This paper introduced a novel framework for designing micro-transit zones, featuring a modified shareability graph model and an efficient two-phase algorithm. Experiments using real-world and synthetic data validated the approach, demonstrating its effectiveness in identifying high-demand zones and its improved performance over a natural baseline method.

One important future direction is to conduct further analysis on the growth of the number of zones as a function of the problem parameters, i.e., provide explicit bounds on the number of zones to consider under various assumptions regarding the network characteristics. Another important future direction is to tackle the problem of optimal micro-transit zones when allowing for inter-zone travel via traditional fixed-route transit, enabling micro-transit to function as both an intra-zone service and the first-and-last-mile connector to mass public transit.

While our experiments demonstrate tractability for medium-sized cities, larger metropolitan areas may require further computational optimizations such as lazy constraint generation, parallelization of clique generation, and memory-efficient data structures. Additionally, real-world implementations should consider metrics beyond maximizing coverage, such as likelihood that each potential traveler uses the service, by extending the objective function to incorporate weighted demand functions that prioritize access in certain regions, ensuring micro-transit serves as a tool for maximizing the social welfare.

\section*{Acknowledgment}
This work was partially supported by the Aizen Climate Scholars Research Fund, US National Science Foundation under award CMMI 2144127 and the US Department of Energy Vehicle Technologies Office under award DE-EE0009212.

\end{document}